\documentclass[11pt]{amsart}
\usepackage{amsmath,amssymb,amsfonts,url}
\numberwithin{equation}{section}

\newtheorem{thm}{Theorem}[section]
\newtheorem{lem}{Lemma}[section]
\newtheorem{rem}{Remark}[section]
\newtheorem{prop}{Proposition}[section]

\begin{document}
\title[Liouville system with singularity]{Classification of Radial Solutions to  Liouville Systems with Singularities} \subjclass{35J60, 35J55}
\keywords{Liouville system, classification of solutions}

\author{Chang-shou Lin}
\address{Taida Institute of Mathematical Sciences\\
and Center for Advanced Study in Theoretical Sciences\\
        National Taiwan University\\
         Taipei 106, Taiwan } \email{cslin@math.ntu.edu.tw}

\author{Lei Zhang}
\address{Department of Mathematics\\
        University of Florida\\
        358 Little Hall P.O.Box 118105\\
        Gainesville FL 32611-8105}
\email{leizhang@ufl.edu}
\thanks{Zhang is supported in part by NSF Grant 0900864 (1027628)}

\date{\today}

\begin{abstract}
Let $A=(a_{ij})_{n\times n}$ be a nonnegative, symmetric, irreducible and invertible matrix. We prove the existence and uniqueness of radial solutions to the following Liouville system with singularity:
$$\left\{\begin{array}{ll}
\Delta u_i+\sum_{j=1}^n a_{ij}|x|^{\beta_j}e^{u_j(x)}=0,\quad \mathbb R^2, \quad i=1,...,n\\
\\
\int_{\mathbb R^2}|x|^{\beta_i}e^{u_i(x)}dx<\infty, \quad i=1,...,n
\end{array}
\right.
$$
where $\beta_1,...,\beta_n$ are constants greater than $-2$. If all $\beta_i$s are negative we prove that all solutions are radial and
 the linearized system is non-degenerate.
\end{abstract}

\maketitle

\section{Introduction}
In this article we consider the following singular Liouville system
\begin{equation}\label{mainsyss}
\left\{\begin{array}{ll}
\Delta u_i+\sum_{j=1}^n a_{ij}|x|^{\beta_j}e^{u_j(x)}=0,\quad \mathbb R^2, \quad i\in I:=\{1,...,n\},\\
\\
\int_{\mathbb R^2}|x|^{\beta_i}e^{u_i(x)}dx<\infty, \quad i\in I.
\end{array}
\right.
\end{equation}
where $\beta_1,..,\beta_n$ are constants greater than $-2$, $A=(a_{ij})_{n\times n}$ is a constant matrix that satisfies
$$(H1): \mbox{ A is symmetric, nonnegative, irreducible and invertible}. $$
$A$ is irreducible means there is no disjoint partition of $I$ into $I_1$ and $I_2$ such that $a_{ij}=0$ for all $i\in I_1$ and $j\in I_2$.
For the system (\ref{mainsyss}), the irreducibility of $A$ means (\ref{mainsyss}) can not be written as two independent subsystems.
If $n=1$ and $a_{11}=1$, the singular Liouville system is reduced to the following single Liouville equation:
\begin{equation}\label{tarane1}
\Delta u+|x|^{\beta}e^{u}=0, \quad \mathbb R^2, \quad \int_{\mathbb R^2} |x|^{\beta} e^u<\infty.
\end{equation}
Prajapat-Tarantello \cite{prajapat} classified all the solutions to (\ref{tarane1}) and proved, on one hand, that if $\beta/2\not \in \mathbb N$, all solutions are radial and can be written as
$$u(x)=\log \frac{\mu}{(1+\frac{\mu}{8(\beta/2+1)^2}|x|^{\beta+2})^2},\quad \mu>0.  $$
On the other hand, a solution may not be symmetric around any point if $\beta/2\in \mathbb N$.
The proof of Prajapat-Tarantello uses properties of integrable system. However, the Liouville system is not integrable and we have to apply new methods.
The purpose of this paper is to prove a classification theorem for all the {\bf radial } solutions to (\ref{mainsyss}).
Let $u=(u_1,...,u_n)$ be a solution to (\ref{mainsyss}) we use $\sigma=(\sigma_1,...,\sigma_n)$ to denote its energy:
\begin{equation}\label{sigma}
\sigma_i=\frac 1{2\pi}\int_0^{\infty}|x|^{\beta_i}e^{u_i(x)}dx,\quad i\in I:=\{1,...,n\}
\end{equation}
and we set $\Lambda_I(\sigma)$ as
$$\Lambda_I(\sigma)=2\sum_{i\in I}(\beta_i+2)\sigma_i-\sum_{i,j\in I} a_{ij}\sigma_i\sigma_j. $$
For $J\subset I$, $\Lambda_J(\sigma)$ is understood similarly.
The main theorem of this article is
\begin{thm}\label{thm1}
Let $A$ satisfy $(H1)$, $\beta_1,...,\beta_n>-2$ be constants,
\begin{enumerate}
\item If $u=(u_1,...,u_n)$ is a radial solution to (\ref{mainsyss}), then
\begin{equation}\label{energyieq}
\Lambda_I(\sigma)=0,\,\, \Lambda_J(\sigma)>0 \,\, \forall \emptyset\subsetneq J\subsetneq I.
\end{equation}
\item
For each $\sigma=(\sigma_1,...,\sigma_n)$ satisfying (\ref{energyieq}), there exists a global radial solution $u$ whose energy is $\sigma$.
\item
If $u$ and $v$ are both radial solutions to (\ref{mainsyss}) with
$$\int_0^{\infty}r^{\beta_i+1}e^{u_i(r)}dr=\int_0^{\infty}r^{\beta_i+1}e^{v_i(r)}dr,\quad i\in I. $$
Then $u_i(r)=v_i(\delta r)+(2+\beta_i)\log \delta$ for some $\delta>0$ and all $i\in I$.
\end{enumerate}
\end{thm}

System (\ref{mainsyss}) is reduced to the following form if $\beta_1=...=\beta_n=0$,
\begin{equation}\label{nosing}
\left\{\begin{array}{ll}
\Delta u_i+\sum_{j\in I}a_{ij}e^{u_j}=0,\quad \mathbb R^2, \\
\\
\int_{\mathbb R^2}e^{u_i}<\infty, \quad \mathbb R^2.
\end{array}
\right.
\end{equation}
Under the assumption $(H1)$ on $A$,
a standard moving-plane argument shows that all $u_1,..,u_n$ are radially symmetric with respect to a common point (see \cite{CSW} for the proof).
The classification of all solutions to (\ref{nosing}) has been completed through the works of Chipot-Shafrir-Wolansky \cite{CSW,CSW1} and the authors \cite{linzhang1}. Among other things Chipot et. al. prove that

\medskip

\emph{Theorem A: (Chipot-Shafrir-Wolansky) Suppose $A$ satisfies $(H1)$, for any solution $u=(u_1,...,u_n)$ to (\ref{nosing}), its energy $\sigma=(\sigma_1,...,\sigma_n)$ belongs to the hypersurface
$$\Gamma:=\{\sigma=(\sigma_1,...,\sigma_n);\,\,  \Lambda_J(\sigma)>0,\,\, \Lambda_J(\sigma)>0, \,\, \forall \emptyset
\subsetneq J\subsetneq I. \,\, \} $$
On the other hand, for any $\sigma\in \Gamma$, there is a solution $u$ of (\ref{nosing}) whose energy is $\sigma$. }

It can be readily verified that the energy of a solution of (\ref{nosing}) is invariant under rigid translations and appropriate scalings: Let $u$ be a global solution to (\ref{nosing}), then $v=(v_1,...,v_n)$ defined by
\begin{equation}\label{11oct4e1}
v_i(y)=u_i(\delta y+x_0)+2\log \delta, \quad i\in I
\end{equation}
for any $x_0\in \mathbb R^2$ and any $\delta>0$ clearly satisfies $\int_{\mathbb R^2}e^{v_i}=\int_{\mathbb R^2}e^{u_i}$ for all $i\in I$. It turns out that for any $\sigma\in \Gamma$, all the global solutions that have the energy $\sigma$ are related by a translation and a scaling described in (\ref{11oct4e1}):

 \medskip

\emph{Theorem B: (\cite{linzhang1}) Suppose $A$ satisfies $(H1)$. Let $u=(u_1,...,u_n)$ and $v=(v_1,...,v_n)$ be global solutions to (\ref{nosing}) such that
$\int_{\mathbb R^2}e^{u_i}=\int_{\mathbb R^2}e^{v_i}$ for all $i\in I$, then $v$ and $u$ are related by (\ref{11oct4e1})
for some $\delta>0$ and $x_0\in \mathbb R^2$. }

\medskip

Theorem A and Theorem B together give a classification of all the solutions to (\ref{nosing}). One obvious question that Theorem \ref{thm1}
raises is, for what $\beta_1,..,\beta_n$ do all the solutions to (\ref{mainsyss}) have to be radially symmetric? We give an affirmative answer for the case of non-positive $\beta$.

\begin{thm}\label{beta0}
Let $u$ be a solution to (\ref{mainsyss}), $A$ satisfy $(H1)$. Suppose $\beta_i\in (-2,0]$ for $i\in I$ and are not all equal to $0$. Then all
components of $u$ are radial functions.
\end{thm}

Systems (\ref{mainsyss}) and (\ref{nosing}) and their reductions appear in many disciplines of mathematics and have profound background in Physics, Chemistry and Ecology.  When (\ref{nosing}) is reduced to one equation, it becomes the classical Liouville equation
$$-\Delta u = e^u, $$
which is related to finding a metric with constant Gauss curvature. In Physics, the Liouville equation represents the electric potential induced by the charge carrier in electrolytes theory \cite{rubinstein} and is closely related to the abelian model in the Chern-Simons theories \cite{jw1,jw2,jost-lin-wang}.

The Liouville systems (\ref{nosing})(\ref{mainsyss}) are used to describe models in the theory of chemotaxis \cite{childress,keller},
in the physics of charged particle beams
\cite{bennet,debye,kiessling2}, and in the theory of
semi-conductors \cite{mock}. For applications of Liouville
systems, see \cite{chanillo2,CSW,linzhang1,linzhang2} and the references therein. Here we note that Liouville systems with singularities are of special importance in Physics and Geometry. For example, the single equation (\ref{tarane1}) appeared in \cite{nolasco2} as a limiting equation in the blow-up analysis of periodic vortices for the Chern-Simons theory of Jackiw and Weinberg \cite{jackiw} and Hong et. al. \cite{hong}. In geometry (\ref{mainsyss}) is related to finding metric with conic singularities \cite{chang1,chang2,chenlin0,kw}.

It is well known that classification theorems are
closely related to blowup analysis and degree-counting theorems. For many equations the asymptotic behavior of blowup solutions are approximated by global solutions. For example, for the Liouville equation
$$\Delta u + V e^{u}=0, \quad \Omega \subset \mathbb R^2, $$
if $V$ is a positive smooth function, blowup solutions near a blowup point can be well approximated by global solutions to
$$\Delta u + e^{u}=0, \quad \mathbb R^2 $$
see \cite{cheng-li-duke1, licmp, chenlin1, zhangcmp}.  If $V$ is nonnegative and the blowup point happens to be a zero of $V$, the profile of blowup solutions is similar to that of the global solutions of (\ref{tarane1}), see \cite{prajapat, bclt,zhangccm}. We expect Theorem \ref{thm1} to be useful in the study of singular Liouville systems defined on Riemann surfaces or domains in $\mathbb R^2$.

The proof of the uniqueness part of Theorem \ref{thm1} (the third statement) is motivated by the authors' previous work \cite{linzhang1} on the Liouville system with no singularity. The existence part (the second statement) is based on the uniqueness result and is therefore significantly different from the duality method used by Chipot. et. al. in \cite{CSW}.
The first statement in Theorem \ref{thm1} is similar to the corresponding case in \cite{CSW}.

For many applications, especially on the construction of bubbling solutions it is important to study the nondegeneracy of the
linearized system. Our next result is concerned with the case when $\beta$ is non-positive.

\begin{thm}\label{nondegen} Let $\beta_i\in (-2,0)$ for all $i\in I$, $u=(u_1,...,u_n)$ solve (\ref{mainsyss}) corresponding to $\beta=(\beta_1,..,\beta_n)$.
Let $\phi=(\phi_1,..,\phi_n)$ be a bounded solution to
$$\Delta \phi_i+\sum_{j\in I} a_{ij} |y|^{\beta_j}e^{u_j(y)}\phi_j(y)=0, \quad \mathbb R^2, i\in I. $$
Then there exists $C\in \mathbb R$ such that $\phi_i(r)=C(ru_i'(r)+2+\beta_i)$ for all $i\in I$.
\end{thm}
\begin{rem} By Theorem \ref{beta0} $u$ is radial in Theorem \ref{nondegen}.
\end{rem}

The organization of this paper is as follows. In section two we list standard tools to be used in the proof of Theorem \ref{thm1}. Then in
section three we prove the three statements of Theorem \ref{thm1}. Theorem \ref{beta0} and Theorem \ref{nondegen} are proved in section four and section five, respectively. Finally in the appendix we provide proofs for the tools used in the proof of Theorem \ref{thm1}.

\medskip

\noindent{\bf Acknowledgement: } Part of the work was finished when the second author was visiting Taida Institute for Mathematics Sciences (TIMS) in March 2011. He would like to thank TIMS for the warm hospitality. The second author is also partially supported by a grant from National Science Foundation.

\section{Preliminary results}

In this section we list a few ODE lemmas to be used in the proof of Theorem \ref{thm1}. Since these lemmas are standard we put their proofs in the appendix, in order not to disturb the main part of the paper.

\begin{lem}\label{global1}
Let $u=(u_1,...,u_n)$ be a solution to (\ref{mainsyss}) where $A$ satisfies $(H1)$. Then
$$u_i(x)=-m_i\log |x|+c_i+o(|x|^{-\delta}),\quad i\in I, \quad |x|>1, $$
$$\nabla u_i(x)=-m_i x/|x|^2+O(|x|^{-\delta-1}), \quad i\in I, \quad |x|>1. $$
where
$$m_i=\sum_{j=1}^n a_{ij}\sigma_j > 2+\beta_i , \quad i\in I, $$
$$c_i=u_i(0)+\int_0^{\infty}\log r\sum_{j=1}^n a_{ij}r^{\beta_j+1}e^{u_j(r)}dr$$
$\delta$ is some positive small number.
\end{lem}

\begin{rem} $u$ is not assumed to be radial in Lemma \ref{global1}.
\end{rem}

\medskip

The next lemma is on the linearized system of (\ref{mainsyss}) expanded along a radial solution $u$:
\begin{equation}\label{lineare1}
(r \phi_i')'+\sum_j a_{ij} r^{\beta_j} e^{u_j} r \phi_j =0, \quad i\in I.
\end{equation}

\begin{lem}\label{lineargrowth}
Let $\phi=(\phi_1,...,\phi_n)$ satisfy (\ref{lineare1}) with $\beta_i>-2$ for all $i\in I$, then $\phi_i(r)=O(\log r)$ at infinity for $i\in I$.
\end{lem}

\begin{lem}\label{ode1}
Let $A$ satisfy $(H1)$, $\beta_i>-2$ for $i\in I$, then for any $c_1,...,c_n\in \mathbb R$, there is a unique solution to
\begin{equation}\label{11sep9e1}
\left\{\begin{array}{ll}
u_i''(r)+\frac 1ru_i'(r)+\sum_{j=1}^n a_{ij}r^{\beta_j}e^{u_j(r)}=0,\quad i=1,..,n,\\
u_i(0)=c_i, \quad i=1,..,n
\end{array}
\right.
\end{equation}
that exists for all $r>0$.
\end{lem}

\begin{rem} $u$ may not have finite energy.
\end{rem}

If we further know that $a_{ii}>0$ for all $i$, then the solution has a finite energy:
\begin{lem}\label{ode2} Let $a_{ii}>0$ for all $i\in I$, then for all $c_1,...,c_n\in \mathbb R$, there exists a solution to
\begin{equation}\label{11sep12e1}
\left\{\begin{array}{ll}
u_i''+\frac 1ru_i'(r)+\sum_{j\in I} a_{ij} r^{\beta_j} e^{u_j(r)}=0, \quad 0<r<\infty, \quad i\in I, \\  \\
\int_0^{\infty} e^{u_i(r)}r^{\beta_j+1}dr<\infty, \quad i\in I, \\  \\
u_i(0)=c_i,\quad i\in I.
\end{array}
\right.
\end{equation}
\end{lem}

\begin{lem}\label{linearlem2}
Let $\phi$ be a solution of
$$\left\{\begin{array}{ll}
(r\phi_i'(r))'+\sum_{j=1}^n a_{ij}r^{\beta_j+1}e^{u_j}\phi_j(r)=0, \quad 0<r<\infty, \\
\\
\phi_i(0)=0,\quad \forall i\in I.
\end{array}
\right.
$$
Suppose $\beta_i>-2$ for all $i\in I$, then
 $\phi_i\equiv 0$ for all $i\in I$.
\end{lem}

\section{Proof of Theorem \ref{thm1}:}

\subsection{ The proof of the first statement of Theorem \ref{thm1}. }

\begin{lem}\label{global2}
Let $u=(u_1,..,u_n)$ be a radial solution of (\ref{mainsyss}) with $A$ satisfying $(H1)$ and $\beta_i>-2$ for all $i$. Then
$$\Lambda_I(\sigma)=0, \quad
\Lambda_J(\sigma)>0,\quad \forall \emptyset \subsetneq J \subsetneq I. $$
\end{lem}

\noindent{\bf Proof of Lemma \ref{global2}:} This proof uses the same idea as in \cite{CSW}.
Let $\tilde u_i(t)=u_i(e^t)$, then
\begin{equation}\label{11oct6e1}
\tilde u_i'(t)\to -m_i\quad \mbox{ as }\quad t\to \infty.
\end{equation}
 The equation for $\tilde u_i(t)$ is
\begin{equation}\label{11sep15e1}
\tilde u_i''(t)+\sum_{j=1}^n a_{ij}e^{(2+\beta_j)t+\tilde u_j(t)}=0, \quad t\in \mathbb R, \quad i\in I.
\end{equation}
Let $z_i(t)=\sum_{j=1}^n a^{ij}\tilde u_j(t)$, then $z_i'(t)\to -\sigma_i$ as $t\to \infty$.
(\ref{11sep15e1}) can be rewritten as
\begin{equation}\label{11sep15e2}
z_i''(t)=-e^{(2+\beta_i)t+\sum_{j=1}^n a_{ij}z_j}, \quad i\in I.
\end{equation}
Clearly $z_i''(t)<0$ for all $i\in I$ and all $t\in \mathbb R$. Let $w_i(t)=z_i'(t)$, then by (\ref{11sep15e2}) and (\ref{11oct6e1})
$$w_i(-\infty)=0, \quad w_i(t)<0 \,\, \forall t, \quad w_i(\infty)=-\sigma_i. $$
In addition we have $w_i'(-\infty)=w_i'(\infty)=0$.
Using the definition of $w_i$ we differentiate (\ref{11sep15e2}) to obtain
\begin{equation}\label{11sep15e3}
w_i''(t)=(2+\beta_i)w_i'(t)+w_i'(t)\sum_{j\in I}a_{ij}w_j(t).
\end{equation}
Taking the summation for $i\in I$ in (\ref{11sep15e3}) we can write (\ref{11sep15e3}) as
$$\sum_{i\in I}w_i''(t)-\sum_{i\in I}(2+\beta_i)w_i'(t)=\sum_{i,j\in I}\frac 12 a_{ij}(w_i(t)w_j(t))'. $$
Integrating $t$ from $-\infty$ to $\infty$ we obtain $\Lambda_I(\sigma)=0$.
For $J\subset I$, summation for $i\in J$ in (\ref{11sep15e3}) leads to
$$\sum_{i\in J}w_i''(t)-\sum_{i\in J}(2+\beta_i)w_i'(t)-\frac 12\sum_{i,j\in J}a_{ij}(w_iw_j)'(t)=\sum_{i\in J,j\in I\setminus J}
a_{ij}w_i'(t)w_j(t). $$
Integrating the above for $t\in (-\infty,\infty)$ we obtain
$$\sum_{i\in J}(2+\beta_i)\sigma_i-\frac 12\sum_{i,j\in J}a_{ij}\sigma_i\sigma_j=\sum_{i\in J, j\in I\setminus J}
a_{ij}\int_{-\infty}^{\infty}w_i'(t)w_j(t)dt. $$
Since the irreducibility of $A$ means that there exists $a_{ij}>0$ for $i\in J, j\in I\setminus J$ we see the right hand side of the above is strictly positive because $w_i'(t)<0$ for all $i\in I$ and $t\in \mathbb R$ and $w_i(t)<0$ for all $i\in I$ and $t\in \mathbb R$. Thus we have obtained $\Lambda_J(\sigma)>0$. Lemma \ref{global2} is established. $\Box$

\medskip

\subsection{ The proof of the third statement of Theorem \ref{thm1}.}
\begin{prop}\label{uniqlinear}
Let $u$ and $v$ both be radial solutions to (\ref{mainsyss}) such that
$$\int_0^{\infty}r^{\beta_i+1}e^{u_i(r)}dr=\int_0^{\infty}r^{\beta_i+1}e^{v_i(r)}dr,\quad i\in I. $$
Then $u_i(r)=v_i(\delta r)+(2+\beta_i)\log \delta$ for some $\delta>0$ and all $i\in I$.
\end{prop}

To prove Proposition \ref{uniqlinear} we first establish a uniqueness result for the linearized system:
\begin{lem}\label{lem1}
Let $\phi=(\phi_1,...,\phi_n)$ be a bounded solution of (\ref{lineare1}), then $\phi_i(r)=C(ru_i'(r)+2+\beta_i)$ for all $i\in I$.
\end{lem}

\noindent{\bf Proof of Lemma \ref{lem1}:} Let
$$\phi^0=(ru_1'(r)+2+\beta_1,...,ru_n'(r)+2+\beta_n). $$
Then by computation $\phi^0$ is a solution to the linearized system. Suppose there exists another bounded solution $\phi^1$ which is not a multiple of $\phi^0$. Without loss of generality we assume $\phi_1^1(0)=0$, as by Lemma \ref{linearlem2} one of $\phi_i^1(0)$ must be different from $2+\beta_i$. To derive a contradiction we set
\begin{eqnarray*}
S=\{\alpha;\quad \exists \mbox{ a bounded solution}\,\, \phi=(\phi_1,...,\phi_n) \mbox{ such that } \phi_1(0)=2+\beta_1, \\
\phi_i(0)=\alpha_i\le 3+\beta_i;,\quad i=2,...,n, \quad \alpha=\min\{2+\beta_1,\alpha_2,...,\alpha_n\}, \\
\int_0^r e^{u_i(s)}\phi_i(s) s^{1+\beta_i} ds>0, \quad \forall r>0, \quad i=1,...,n. \quad \}.
\end{eqnarray*}
First we see that $2+\min\{\beta_1,...,\beta_n\}\in S$. Indeed the expression of $\phi^0$ gives
$$\int_0^r s^{1+\beta_i}e^{u_i(s)}\phi_i^0(s)ds=r^{2+\beta_i}e^{u_i(r)}>0. $$
Next we observe that $S$ is a bounded set. Indeed, suppose $\alpha<0$ is in $S$, let $\tilde \phi$ be the function corresponding to $\alpha$, then $\exists j\in I$ such that $\tilde \phi_j(0)=\alpha$. This leads to $\int_0^r s^{1+\beta_j}e^{u_j(s)}\tilde \phi_j(s)ds<0$ for $r$ small,
a contradiction to the definition of $S$.
Let $\bar \alpha$ be the infimum of $S$ and let $\alpha^k=(\alpha_1^k,...,\alpha_n^k)\in S$ be a sequence in $S$ that tends to $\bar \alpha$ from above. Suppose $\phi^k=(\phi_1^k,...,\phi_n^k)$ is the solution corresponding to $\alpha^k$, then we claim  that $\phi^k$ converges to $\bar \phi=(\bar \phi_1,...,\bar \phi_n)$, which is also a bounded solution with the strict monotonicity property described in $S$. Indeed, let $\psi^m=(\psi_1^m,...,\psi_n^m)$ be the solution to the linearized system such that $\psi_i^m(0)=\delta_i^m$. Then by Lemma \ref{linearlem2}
$$\phi^k=\sum_{m=1}^n \alpha_m^k\psi^m. $$
Here we recall that by Lemma \ref{lineargrowth} $\psi_i^m(r)=O(\log r)$ for $r$ large.
Since $\bar \alpha\le \alpha_i^k\le 3+\beta_i$ for $i\in I$ and all $k$. Along a subsequence $\phi^k$ tends to $\bar \phi$ over all compact subsets of $\mathbb R$. The monotonicity property of $\phi^k$ implies
$$\int_0^r e^{u_i(s)}\bar \phi_i(s) s^{1+\beta_i} ds\ge 0, \quad \forall i\in I, \quad \forall r>0. $$

On the other hand, since $\phi^k$ are all bounded functions, for each $\phi_i^k$ we find $r_l\to \infty$ as $l\to \infty$ such that
$r_l(\phi_i^k)'(r_l)\to 0$. From the equation for $r\phi_i^k$ we have
$$\int_0^{\infty}\sum_{j=1}^n a_{ij}r^{\beta_j+1}e^{u_j(r)}\phi_j^k(r)dr=0,\quad \forall i\in I. $$
Since $A$ is invertible
$$0=\int_0^{\infty}e^{u_i(r)}\phi_i^k(r)r^{\beta_i+1}dr=\sum_{m=1}^n\alpha_m^k\int_0^{\infty}e^{u_i(r)}\psi_i^m(r)r^{\beta_i+1}dr. $$
Since $\psi_i^m(r)=O(\log r)$, $\displaystyle{\int_0^{\infty}e^{u_i(r)}\psi_i^m(r)r^{\beta_i+1}dr}$ is well defined, we
let $\alpha^k\to (\bar \alpha_1,...,\bar \alpha_n)$ to obtain
\begin{equation}\label{11sep8e1}
\int_0^{\infty}e^{u_i(s)}\bar \phi_i(s)s^{\beta_i+1}ds=0,\quad \forall i\in I.
\end{equation}
As a consequence of (\ref{11sep8e1}), $\bar \phi$ is bounded. Indeed, the equation for $\bar \phi$ is
$$ (r\bar \phi_i'(r))'=-\sum_j a_{ij} r^{\beta_j+1} e^{u_j(r)} \bar \phi_j(r), \quad r>0. $$
Using $\bar \phi_i(r)=O(\log r)$, $r^{\beta_i+2}e^{u_i(r)}=O(r^{-\delta})$ for some $\delta>0$ (Lemma \ref{global1}) and (\ref{11sep8e1}) we know
$$\int_0^r e^{u_i(s)}\bar \phi_i(s)s^{\beta_i+1}ds=0-\int_r^{\infty}e^{u_i(s)}\bar \phi_i(s)s^{\beta_i+1}ds =O(r^{-\delta/2}) $$
for $r$ large. Thus $\bar \phi_i'(r)=O(r^{-1-\delta})$ for all $r$ large, which implies that $\bar \phi_i$ is bounded.
Since each $\bar \phi_i$ is a non-increasing function, (\ref{11sep8e1}) implies that $\bar \phi_i$ decreases to a negative constant when $r\to \infty$. Indeed, by (\ref{11sep8e1}) either $\bar \phi_i\equiv 0$ or $\bar \phi_i$ decreases to a negative constant. The first possibility does not exist, because the fact $\bar \phi_1(0)=2+\beta_1>0$ implies that $\bar \phi_1$ decreases into a negative constant at infinity. Also
$\int_0^r s^{1+\beta_1}e^{u_1(s)}\bar \phi_i(s)ds>0$ for all $r$. Consequently for all $i$ in the set $I_1:=\{i\in I;\,\, a_{i1}>0\, \}$,
$$ r\bar \phi_i'(r)\le -a_{i1}\int_0^r s^{1+\beta_1}e^{u_1(s)}\bar \phi_1(s)ds<0, \quad \forall r>0. $$
Therefore $\bar \phi_i$ strictly decreases to a negative constant for all $i\in I_1$. We can further define
$$I_2:=\{i\in I;\quad a_{ij}>0\,\, \mbox{ for some } j\in I_1. \,\, \}. $$
By the same reason as above $\bar \phi_i$ decreases to a negative constant at infinity for all $i\in I_2$. By the irreducibility of $A$ all the components of $\bar \phi$ decrease to negative constants at infinity.

Now we claim that $\bar \alpha-\epsilon\in S$ for $\epsilon>0$ small. To see this, consider $\bar \phi+t \phi^1$ for $|t|$ sufficiently small. Recall that $\phi_1^1(0)=0$, thus $\bar \phi_1(0)+t\phi_1^1(0)=2+\beta_1$. Clearly $\bar \phi+t\phi^1$ solves (\ref{lineare1}). By choosing $t$ positive or negative with $|t|$ small we can make
$$\min_{i\in I} \bar \phi_i(0)+t\phi_i^1(0)=\bar \alpha-\epsilon>0. $$
Since $\bar \phi+t \phi^1$ is bounded we have
$$\int_0^{\infty}e^{u_i}(\bar \phi_i+t \phi_i^1)s^{\beta_i +1}ds=0,\quad i=1,...,n. $$
Since $\bar \phi_i(r)$ tends to a negative constant as $r\to \infty$ and $\phi^1$ is bounded, we know for $r$ large and $|t|$ small
$$\int_r^{\infty}e^{u_i}(\bar \phi_i(s)+t\phi_i^1(s))s^{\beta_i+1}ds<0. $$
Consequently
$$\int_0^r e^{u_i(s)}(\bar \phi_i(s)+t\phi_i^1(s))s^{\beta_i+1}ds>0 \quad \forall r>0. $$
Thus $\bar \alpha-\epsilon\in S$ for some $\epsilon>0$ small, a contradiction to the definition of $\bar \alpha$. Lemma \ref{lem1} is established. $\Box$

\medskip

\noindent{\bf Proof of Proposition \ref{uniqlinear}:}  We shall consider
\begin{equation}\label{11sep17e1}
\left\{\begin{array}{ll}
u_i''(r)+\frac 1r u_i'(r)+\sum_{j\in I} a_{ij} r^{\beta_j}e^{u_j(r)}=0, \quad 0<r<\infty, \\
\\
\int_0^{\infty} e^{u_i(r)}r^{\beta_i+1}dr<\infty, \quad \forall i\in I, \\
\\
u_i(0)=c_i,\quad i=1,...,n-1, \quad u_n(0)=0.
\end{array}
\right.
\end{equation}

Let
$$\Pi_2:=\{\sigma=(\sigma_1,..,\sigma_n);\quad \Lambda_I(\sigma)=0,\quad \Lambda_J(\sigma)>0, \quad \forall \emptyset\subsetneq J\subsetneq I. \quad \}. $$
$$\Pi_1:=\{{\bf C}=(c_1,..,c_{n-1});\quad (\ref{11sep17e1}) \mbox{ has a solution}. \quad \}. $$
Note that by Lemma \ref{ode2} $\Pi_1=\mathbb R^{n-1}$ if $a_{ii}>0$ for all $i$. We claim that the mapping from $\Pi_1$ to $\Pi_2$ is locally one to one. Indeed, let ${\bf M}$ be the following matrix:
$${\bf M}=\left(\begin{array}{ccc}
\partial_{c_1}\sigma_1 & ... & \partial_{c_{n-1}}\sigma_1 \\
...  &  ... & ... \\
\partial_{c_1} \sigma_{n-1}  & ...  &   \partial_{c_{n-1}}\sigma_{n-1}
\end{array}
\right )
$$
We claim that ${\bf M}$ is nonsingular. We prove this claim by contradiction. Suppose there exists a non-zero vector ${\bf D}=(d_1,...,d_{n-1})^T$ such that ${\bf MD}=0$. Then by setting $\gamma=d_1c_1+...+d_{n-1}c_{n-1}$ we have
\begin{equation}\label{11sep12e2}
\partial_{\gamma}\sigma_1=\partial_{\gamma}\sigma_2=...=\partial_{\gamma}\sigma_{n-1}=0.
\end{equation}
For $\Pi_2$, $\Lambda_I(\sigma)=0$ reads
$$\sum_{i,j\in I}a_{ij}\sigma_i\sigma_j=2\sum_{i\in I}(2+\beta_i)\sigma_i. $$
By differentiating both sides with respect to $\gamma$ we have
$$\sum_i (\sum_j a_{ij}\sigma_j-2-\beta_i)\partial_{\gamma}\sigma_i=0. $$
Since $\Lambda_J(\sigma)>0$ implies $\sum_j a_{ij}\sigma_j>2+\beta_i$, (\ref{11sep12e2}) implies $\partial_{\gamma}\sigma_n=0$. Set $\phi_i=\partial_{\gamma}u_i$ ($i\in I$), then $\phi=(\phi_1,...,\phi_n)$ satisfies (\ref{lineare1}) and
$$
\phi_i(0)=d_i,\quad i=1,...,n-1,\quad \phi_n(0)=0. $$ From $\partial_{\gamma}\sigma_i=0$ ($i\in I$) we have
\begin{equation}\label{11oct6e8}
\int_0^{\infty}e^{u_i}\phi_i(s)s^{1+\beta_i} ds=0,\quad i\in I.
\end{equation}
As a consequence of (\ref{11oct6e8}), $\phi$ is bounded. Indeed, integrating (\ref{lineare1}) from $0$ to $r$
\begin{eqnarray*}
r\phi_i'(r)=-\int_0^r \sum_j a_{ij}s^{1+\beta_j}e^{u_j(s)}\phi_j(s)ds\\
=\int_r^{\infty} a_{ij}s^{1+\beta_j}e^{u_j(s)}\phi_j(s)ds
=O(r^{-\delta})
\end{eqnarray*}
for some $\delta>0$. Therefore $\phi'(r)=O(r^{-1-\delta})$, which proves that $\phi_i$ is bounded.
By Lemma \ref{lem1} $\phi_i=c(ru_i'+2+\beta)$, then we see immediately that $c=0$ because $\phi_n(0)=0$, this is not possible because not all $d_i$'s are zero.  Therefore we have proved that ${\bf M}$ is nonsingular for all
${\bf C}=(c_1,...,c_{n-1})\in \Pi_1$.

We further assert that there is one-to-one correspondence between $\Pi_1$ and $\Pi_2$. This is proved in two steps as follows.

Case 1: $a_{ii}>0$, $i\in I$.

In this case, by Lemma \ref{ode2} $\Pi_1=\mathbb R^{n-1}$. The mapping from $\Pi_1$ to $\Pi_2$ is proper and locally one to one. Here we claim that $\Pi_2$ is simply connected. Assuming this, since $\mathbb R^{n-1}$ and $\Pi_2$ are simply connected, there is one to one correspondence between them. Let $u=(u_1,...,u_n)$ and $v=(v_1,..,v_n)$ be two radial solutions such that $u_n(0)=v_n(0)=0$, $\int_{\mathbb R^2}|x|^{\beta_i}e^{u_i}=\int_{\mathbb R^2}|x|^{\beta_i}e^{v_i}$ ($i\in I$). Then $u_i(0)=v_i(0)$
for $i=1,...,n-1$. By Lemma \ref{ode1} $u_i\equiv v_i$ for all $i\in I$. Now we prove that $\Pi_2$ is simply connected. Indeed,
using $m_i=\sum_j a_{ij}\sigma_j$, $\Lambda_I(\sigma)=0$ can be written as
\begin{equation}\label{11oct25e1}
\sum_{i,j\in I} a^{ij}(2+\beta_i)(2+\beta_j)=\sum_{i,j\in I}a^{ij}(m_i-2-\beta_i)(m_j-2-\beta_j).
\end{equation}
Therefore $\Pi_2$ is part of a quadratic surface, the boundary of which is restricted by $\Lambda_{J_i}(\sigma)=0$ where $J_i$ is $I$ with the index $i$ removed. $\Lambda_{J_i}(\sigma)>0$ reads
$$m_i-2-\beta_i>\frac{a_{ii}}2 \sigma_i. $$
In another word in the coordinate system represented by $m_i$, we use $n$ coordinate planes to bound the quadratic hypersurface described in
(\ref{11oct25e1}). Other restrictions $\Lambda_J>0$, when $J$ is obtained from $I$ with at least two indices removed, do not affect the topological information of $\Lambda_I(\sigma)=0$.
Thus $\Pi_2$ is a part of the quadratic hyper-surface in the first quadrant and is therefore simply connected.
Proposition \ref{uniqlinear} is proved in this case.

Case 2: There exists $i_0$ such that $a_{i_0,i_0}=0$. We prove this case by a contradiction. Suppose $c^k=(c_1^k,...,c_{n-1}^k)$ ($k=1,2$) are two distinct points on  $\Pi_1$ that correspond to the same energy: let $u^1,u^2$ be two solutions corresponding to $c^1$ and $c^2$ respectively such that
$$\int_0^{\infty} e^{u_i^1(r)}r^{1+\beta_i}dr=\int_0^{\infty} e^{u_i^2(r)}r^{1+\beta_i}dr=\sigma_i,\quad i\in I. $$
Since the matrix $\displaystyle{\bigg (\frac{\partial \sigma}{\partial c}\bigg )\in \mathbb M_{(n-1)\times (n-1)}}$ is nonsingular at $c^1$ and $c^2$,
there is a one-to-one mapping between a neighborhood of $c^k$ to a neighborhood of $\sigma$ in $\Pi_2$. Since $c^1\neq c^2$, we choose the neighborhoods around them to be disjoint.

Now consider a perturbation system
\begin{equation}\label{per-sys}
\left\{\begin{array}{ll}
u_i''(r)+\frac 1r u_i'(r)+\sum_{j\in I}(a_{ij}+\epsilon \delta_{ij})r^{\beta_j}e^{u_j}=0, \quad r>0, \quad i\in I,\\
\\
\int_0^{\infty} r^{\beta_i+1}e^{u_i}dr<\infty, \quad i\in I, \\
\\
u_1(0)=c_1,...u_{n-1}(0)=c_{n-1},\quad u_n(0)=0.
\end{array}
\right.
\end{equation}
Let $u^{k,\epsilon}$ be the solution to (\ref{per-sys}) that corresponds to the initial condition $c^k=(c_1^k,...,c_{n-1}^k,0)$ ($k=1,2$). Let
$\sigma^{k,\epsilon}=(\sigma^{k,\epsilon}_1,..,\sigma^{k,\epsilon}_n)$ be defined as
$\sigma^{k,\epsilon}_i=\int_0^{\infty}r^{\beta_i+1}e^{u^{k,\epsilon}_i(r)}dr$ ($i=1,..,n$).
We claim that
\begin{equation}\label{11sep13e3}
\sigma^{k,\epsilon}=(\sigma_1,..,\sigma_n)+\circ(1),\quad k=1,2.
\end{equation}
and
\begin{equation}\label{11sep13e4}
\frac{\partial \sigma^{k,\epsilon}_i}{\partial c_j}=\frac{\partial \sigma_i}{\partial c_j}+\circ(1),
\quad i=1,..,n,\quad j=1,..,n-1,\quad k=1,2.
\end{equation}
Assuming (\ref{11sep13e3}) and (\ref{11sep13e4}) for the moment.
Now the matrix
$$\left(\begin{array}{ccc}
\partial_{c_1}\sigma^{k,\epsilon}_1 & \ldots & \partial_{c_{n-1}}\sigma^{k,\epsilon}_1
\\
\vdots &\vdots & \vdots \\
\partial_{c_1}\sigma^{k,\epsilon}_{n-1} & \ldots
&\partial_{c_{n-1}}\sigma^{k,\epsilon}_{n-1}
\end{array}
\right )
$$
is non-singular at $c^k$ ($k=1,2$) for $\epsilon$ small.
On the other hand, $\sigma^{1,\epsilon}$ and $\sigma^{2,\epsilon}$ both satisfy
\begin{equation}\label{feb25e5}
\left\{\begin{array}{ll}
\Lambda^{\epsilon}_I(\sigma^{k,\epsilon}):=\sum_{i\in I}2(2+\beta_i)\sigma^{k,\epsilon}_i-\sum_{i,j\in I}(a_{ij}+\epsilon \delta_{ij})
\sigma^{k,\epsilon}_i\sigma^{k,\epsilon}_j=0 \\
\\
\Lambda^{\epsilon}_J>0,\quad 0\varsubsetneqq  J\varsubsetneqq  I.
\end{array}
\right.
\end{equation}
We use $\Pi^{\epsilon}$ to represent the hyper-surface described as above.
For
$\sigma^{2,\epsilon}=(\sigma^{2,\epsilon}_1,..,\sigma^{2,\epsilon}_n)\in \Pi^{\epsilon}$, we can find
$c^{1,\epsilon}=(c^{1,\epsilon}_1,..,c^{1,\epsilon}_{n-1})$ such that
$$c^{1,\epsilon}_j=c^1_j+\circ(1),\quad j=1,2,..,n-1$$ and a solution $\bar u^{1,\epsilon}$ of (\ref{per-sys})
with the initial condition $(c^{1,\epsilon}_1,..,c^{1,\epsilon}_{n-1},0)$
such that
$$\int_0^{\infty}r^{\beta_j+1}e^{\bar u^{1,\epsilon}_j}dr=\sigma^{2,\epsilon}_j,\quad j=1,2,..,n-1. $$
After using $\Lambda_I^{\epsilon}(\sigma^{2,\epsilon})=0$ in (\ref{feb25e5}) we have
$$\int_0^{\infty}r^{\beta_n+1}e^{\bar u^{1,\epsilon}_n}dr=\sigma^{2,\epsilon}_n. $$
Then the difference
between $c^1$ and $c^2$ implies $c^{1,\epsilon}\neq c^2$ for $\epsilon$ small.
A contradiction to the uniqueness property satisfied by the system (\ref{per-sys}).

To finish the proof we now verify (\ref{11sep13e3}) and (\ref{11sep13e4}).
Here we require $\epsilon\in (0,\delta_0)$ where $\delta_0$ is so small that the matrix $(a_{ij}+\epsilon\delta_{ij})_{n\times n}$
is non-singular for all $\epsilon\in (0,\delta_0)$.

For $u^k$, there exists $R_0$ large such that for $r>R_0$ and some $\delta>0$,
$$(u_i^k)'(r)r\le -2-\beta_i-2\delta, \quad i=1,...,n,\quad k=1,2. $$
For $\delta_0$ small we have $u_i^{k,\epsilon}$ converges uniformly to $u_i^k$ over $0\le r\le R_0$. For $r=R_0$ we have
$$
(u^{k,\epsilon}_j(r))'r\leq -(2+\beta_j+\delta)\quad\mbox{at}\quad
r=R_0,\quad 0\le \epsilon \le \delta_0.
$$
Then by the super-harmonicity of $u^{k,\epsilon}_j$ it is easy to show
$$
(u^{k,\epsilon}_j(r))'r\le -(2+\beta_j+\delta)\quad\mbox{for}\quad r\geq
R_0.
$$

\noindent Thus, $\exists C>0$ and $R_1\geq R_0$ such that
\begin{equation}\label{feb25e3}
r^{\beta_j}e^{u^{k,\epsilon}_j(r)}\leq Cr^{-(2+\delta)}\quad\mbox{for}\quad
r\geq R_1
\end{equation}

Hence for $k=1,2$,
$$
{\sigma_j}^{\epsilon}=\int_0^\infty
e^{u^{k,\epsilon}_j(r)}r^{\beta_j+1}dr=\int_0^\infty
e^{u_j^k(r)}r^{\beta_j+1}dr+o(1)=\sigma_j+\circ(1), \quad j=1,..,n.
$$
(\ref{11sep13e3}) is verified.
To show (\ref{11sep13e4})

\begin{equation}\label{feb25e1}
\frac{\partial \sigma^{\epsilon}_i}{\partial c_j}=\int_0^{\infty}r^{\beta_i+1}e^{u^{k,\epsilon}_i(r)}
\frac{\partial u^{k,\epsilon}_i}{\partial c_j}(r)dr,\quad i=1,..,n,\quad k=1,2.
\end{equation}
$\frac{\partial u^{k,\epsilon}}{\partial c}$ satisfies the
following linearized equation:
$$-\Delta (\frac{\partial u^{k,\epsilon}_i}{\partial c_l})=\sum_{j=1}^2(a_{ij}+\epsilon \delta_{ij})r^{\beta_j}e^{u^{k,\epsilon}_j}
\frac{\partial u^{k,\epsilon}_j}{\partial c_l},\quad i=1,..,n, \quad l=1,...,n-1.
$$
By Lemma \ref{lineargrowth}
\begin{equation}\label{feb25e2}
|\frac{\partial u^{k,\epsilon}_i}{\partial c_l}(r)|\le C\ln r,\quad r\ge 2, \quad i=1,..,n,\quad l=1,..,n-1
\end{equation}
where the constant $C$ is independent of $\epsilon\in (0,\delta_0)$. Moreover, for any fixed $R>0$,
$\frac{\partial u^{1,\epsilon}_i}{\partial c_l}(r)$ converges uniformly to
$\frac{\partial u^1_i}{\partial c_l}(r)$ over $0<r<R$ with respect to $\epsilon$.
Using the decay estimates (\ref{feb25e3}) and (\ref{feb25e2}) in (\ref{feb25e1}) we obtain (\ref{11sep13e4})
by elementary analysis. Proposition \ref{uniqlinear} is proved in all cases. $\Box$

\medskip

\subsection{The proof of the second statement of Theorem \ref{thm1}.}

Our proof is based on the uniqueness result and is completely different from the method employed in \cite{CSW}.
We divide the proof into two cases according to the diagonal entries of $A$.

\noindent{Case one: $a_{ii}>0$ for all $i\in I$. }

In this case, by Lemma \ref{ode2}, for any $c_1,...,c_{n-1}\in \mathbb R$, there exists a unique finite energy solution $u=(u_1,..u_n)$ such that
$u_i(0)=c_i$ for $i=1,..,n-1$ and $u_n(0)=0$. By Proposition \ref{uniqlinear} there is a bijection between the initial condition $(c_1,..,c_{n-1},0)$ and $\Pi_2$ (see the notation in the proof of Proposition \ref{uniqlinear}). Thus Theorem \ref{thm1} is proved in this case.

\medskip

\noindent{ Case two: There exists $i_0\in I$ such that $a_{i_0,i_0}=0$. }

Let $\sigma\in \Pi_2$, then for $\epsilon>0$ we consider
$$\left\{\begin{array}{ll}
-\Delta u_i^{\epsilon}=\sum_{j\in I} (a_{ij}+\epsilon \delta_{ij})|x|^{\beta_j}e^{u_j^{\epsilon}(x)},\quad \mathbb R^2, \\
\\
u_i^{\epsilon}(0)=c_i^{\epsilon},\quad i=1,...,n-1, \quad u_n^{\epsilon}(0)=0, \\
\\
\int_0^{\infty} r^{\beta_i+1}e^{u_i^{\epsilon}(r)}dr=\sigma_i^{\epsilon}, \quad i\in I
\end{array}
\right.
$$
where $\sigma^{\epsilon}=(\sigma_1^{\epsilon},...,\sigma_n^{\epsilon})$ is a point on the hyper-surface
$$
\Pi_2^{\epsilon}:=\{\sigma=(\sigma_1,...,\sigma_n);\,\,  \sigma_i>0, \,\,\, \forall i\in I,\,\,
\Lambda_I^{\epsilon}(\sigma)=0,\,\, \Lambda_J^{\epsilon}(\sigma)>0,\, \forall \emptyset\subsetneq J\subsetneq I\}$$
such that $\sigma^{\epsilon}\to \sigma$ as $\epsilon\to 0$.
Here we recall that $\Lambda_I^{\epsilon}(\sigma)$ is defined as
$$\Lambda_I^{\epsilon}(\sigma):=2\sum_{i\in I}(2+\beta_i)\sigma_i-\sum_{i,j\in I}(a_{ij}+\epsilon \delta_{ij})\sigma_i\sigma_j. $$

The vector $(c_1^{\epsilon},...,c_{n-1}^{\epsilon},0)$ is the initial condition corresponding to $\sigma^{\epsilon}$.
Now we claim that
\begin{equation}\label{11sep23e1}
\max c_i^{\epsilon}\le C \quad i=1,...,n-1
\end{equation}
for some $C>0$ independent of $\epsilon$. Indeed, if this is not the case, without loss of generality we assume $c_1^{\epsilon}$ is the largest
among $c_i^{\epsilon}$ and tends to infinity. Re-scale $u^{\epsilon}$ according to $c_1^{\epsilon}$ to make the maximum of all components at $0$ equal to $0$. The re-scaled system has to converge in $C_{loc}^2(\mathbb R^2)$ norm to a partial system. Indeed, the first component converges because all the components are bounded. The $n-th$ component tends to $-\infty$ because the initial condition is $0$ before the scaling and all components are non-increasing. Therefore for the limit function $v=(v_1,...,v_n)$ without loss of generality we assume $v_{m+1}=...v_n=0$ for some $1<m<n$.
For $i=1,...,m$ we easily observe that
\begin{equation}\label{11sep23e2}
\bar \sigma_i:=\int_0^{\infty}r^{\beta_i+1}e^{v_i(r)}dr\le \sigma_i, \quad i=1,...,m.
\end{equation}
The reason is for each fixed $R>0$ we have
$$ \int_0^{R}r^{\beta_i+1}e^{v_i(r)}dr\le \sigma_i^{\epsilon}+\circ(1), \quad i=1,...,m. $$
Clearly $(v_1,...,v_m)$ satisfies
$$\left\{\begin{array}{ll}
\Delta v_i+\sum_{j=1}^m a_{ij} r^{\beta_j}e^{v_j}=0, \quad i=1,...,m, \\
\\
\int_0^{\infty}r^{\beta_i+1}e^{v_i(r)}dr\le \sigma_i,\quad i=1,..,m.
\end{array}
\right.
$$
By Lemma \ref{global1}
\begin{equation}\label{11sep23e3}
\sum_{j=1}^m a_{ij}\bar \sigma_j>2+\beta_i,\quad i=1,...m.
\end{equation}
We claim that $\bar \sigma=(\bar \sigma_1,...,\sigma_n)$ with $\bar \sigma_{m+1}=...=\bar \sigma_n=0$ satisfies
$\Lambda_I(\bar \sigma)=0$. Indeed, let $v_{m+1}=...=v_n\equiv 0$ and $H_i=1$ if $i=1,...,m$ and $H_i=0$ for $i=m+1,...,n$. Then
the system for $v$ can be written as
$$\Delta v_i+\sum_{j=1}^n a_{ij} r^{\beta_j} H_j e^{v_j}=0, \quad i=1,...,n. $$
Apply the standard method to obtain the Pohozaev identity to the system above we have $\Lambda_I(\bar \sigma)=0$. Let $J=\{1,...m\}$ we have
$\Lambda_J(\bar \sigma)=0$. Let $z_i=\sigma_i-\bar \sigma_i$. From the definition of $\bar \sigma_i$ we know that $z_i\ge 0$ for $i=1,...,m$ Since $1<m<n$ we have $\Lambda_J(\sigma)>0$, $\Lambda_J(\sigma)-\Lambda_J(\bar \sigma)>0$ gives
\begin{equation}\label{11sep23e4}
\sum_{i\in J} \bigg ( (\sum_{j\in J}a_{ij}\sigma_j-(2+\beta_i))z_i +
(\sum_{j\in J}a_{ij}\bar \sigma_j-(2+\beta_i))z_i\bigg )<0.
\end{equation}
Since $\sum_{j\in J}a_{ij}\bar \sigma_j>2+\beta_i$ for all $i\in J$, we also have $\sum_{j\in J}a_{ij}\sigma_j>2+\beta_i$ for all $i\in J$
because $\sigma_i\ge \bar \sigma_i$. Clearly (\ref{11sep23e4}) is impossible. (\ref{11sep23e1}) is proved. Similarly there is a lower bound for
$c_1^{\epsilon},...,c_{n-1}^{\epsilon}$. As $\epsilon\to 0$, the $u^{\epsilon}$ converges to $u$ that corresponds to $\sigma$.
Theorem \ref{thm1} is proved in both cases. $\Box$

\section{Proof of Theorem \ref{beta0}}

The proof of Proposition 4.1 in \cite{CSW} can be readily applied to prove Theorem \ref{beta0}. We include it for the convenience of readers.

For $\lambda>0$, let $u_i^{\lambda}(x_1,x_2)=u_i(2\lambda-x_1,x_2)$. Set $\Sigma_{\lambda}=\{x\in \mathbb R^2;\quad x_1>\lambda. \}$ and $T_{\lambda}$ be the boundary of $\Sigma_{\lambda}$.  The equation for $u^{\lambda}
=(u_1^{\lambda},...,u_n^{\lambda})$ is
\begin{equation}\label{ulambda}
\Delta u_i^{\lambda}+\sum_{j\in I} a_{ij}|x^{\lambda}|e^{u_j^{\lambda}}=0, \quad i\in I
\end{equation}
where $x^{\lambda}=(2\lambda-x_1,x_2)$.
Set $w_i^{\lambda}=u_i^{\lambda}-u_i$ to be defined in $\Sigma_{\lambda}$ for $\lambda>0$.
For $w_i^{\lambda}$ we have
$$\Delta w_i^{\lambda}+\sum_j a_{ij}|x|^{\beta_i}e^{\xi_j^{\lambda}}w_j^{\lambda}=-\sum_j a_{ij}(|x^{\lambda}|^{\beta_j}-|x|^{\beta_j})e^{u_j^{\lambda}} $$
where
$$e^{\xi_i^{\lambda}}=\frac{e^{w_i^{\lambda}}-e^{w_i}}{w_i^{\lambda}-w_i}=\int_0^1e^{w_i+t(w_i^{\lambda}-w_i)}dt. $$
Since $\beta_i\le 0$ for all $i\in I$,
\begin{equation}\label{11nov11e1}
\Delta w_i^{\lambda}+\sum_j a_{ij}|x|^{\beta_i}e^{\xi_j^{\lambda}}w_j^{\lambda}\le 0.
\end{equation}
Let $f=\log \log (|x|+3)$, then
$$\Delta f(x)=\frac{3}{r(r+3)^2\log (r+3)}-\frac{1}{(r+3)^2\log^2(r+3)}. $$
Therefore for any $\epsilon>0$, there exists $C(\epsilon)>0$ such that
\begin{equation}\label{11nov11e2}
\frac{\Delta f}{f}\le -\frac{1}{r^{2+\epsilon}},\quad r>C(\epsilon).
\end{equation}
Let $z_i^{\lambda}=w_i^{\lambda}/f$, then the following lemma holds.
\begin{lem}\label{lesslem1}
There exists $R>0$ independent of $\lambda$ such that for $\lambda>0$, if $x_0$ is a point where a negative minimum of $\min\{z_1^{\lambda},..
z_n^{\lambda}\}$ is attained, then $x_0\in B_R$.
\end{lem}

\noindent{\bf Proof of Lemma \ref{lesslem1}:}
From (\ref{11nov11e1}) we obtain
\begin{equation}\label{11nov11e3}
\Delta z_i^{\lambda}+2\nabla z_i^{\lambda}\frac{\nabla f}{f}+z_i^{\lambda}\frac{\Delta f}{f}+\sum_j a_{ij}e^{\xi_j^{\lambda}}
z_j^{\lambda}\le 0.
\end{equation}
Suppose $z_i^{\lambda}(x_0)=\min_j z_j^{\lambda}(x_0)<0$ and $x_0$ is where the negative minimum for $z_i^{\lambda}$ is attained. Here we note that the global minimum of $z_i^{\lambda}$ should be attained. Indeed, by Lemma \ref{global1},
$$u_i(x)=-m_i\log |x|+c_i+O(|x|^{-\delta}) $$
when $|x|$ is large. Thus, for $\lambda>0$, since $|x^{\lambda}|<|x|$,
$$w_i^{\lambda}(x)=u_i^{\lambda}(x)-u_i(x)\ge O(|x|^{-\delta}),\quad |x|>>1. $$
Thus $\lim_{|x|\to \infty}z_i^{\lambda}(x)\ge 0$.
Let $J=\{j\in I;
\quad z_j(x_0)\le 0 \}$. Here we observe that the image of the origin is not in $J$, because of the decay rate of $u_i$. We rewrite (\ref{11nov11e3}) as
\begin{equation}\label{11nov11e4}
\Delta z_i^{\lambda}+2\nabla z_i^{\lambda}\frac{\nabla f}{f}+z_i^{\lambda}\frac{\Delta f}{f}+\sum_{j\in J} a_{ij}e^{\xi_j^{\lambda}}
z_j^{\lambda}\le 0.
\end{equation}
in a small neighborhood of $x_0$. Then at $x_0$,
$$\Delta z_i^{\lambda}(x_0)\ge 0,\quad 2\nabla z_i^{\lambda}(x_0)\frac{\nabla f(x_0)}{f(x_0)}
=0. $$ For $j\in J$, since $w_j^{\lambda}(x_0)\le 0$, we have $u_j^{\lambda}(x_0)\le u_j(x_0)$, so if $|x_0|$ is large, by Lemma \ref{global1},
$e^{\xi_j^{\lambda}(x)}\sim |x|^{-2-\delta}$ for $x$ close to $x_0$ and some $\delta>0$. Thus
$$
 \sum_{j\in J}a_{ij}e^{\xi_j^{\lambda}(x_0)}z_j^{\lambda}(x_0)
 \le z_i(x_0)\sum_ja_{ij}e^{\xi_j^{\lambda}(x_0)}\le Cz_i^{\lambda}(x_0)|x_0|^{-2-\delta}.
$$
On the other hand if $|x_0|$ is large
$$z_i^{\lambda}(x_0)\frac{\Delta f(x_0)}{f(x_0)}> |z_i^{\lambda}(x_0)||x_0|^{-2-\epsilon}. $$
Therefore by choosing $\epsilon<\delta/2$ we see that (\ref{11nov11e4}) can not hold if $|x_0|$ is large.
Lemma \ref{lesslem1} is established. $\Box$

\medskip

By Lemma \ref{lesslem1} and Lemma \ref{global1}, $\min\{w_1^{\lambda},...,w_n^{\lambda}\}>0$ in $\Sigma_{\lambda}$ for $\lambda$ sufficiently large. Thus set
$$\bar \lambda:=\inf \{\lambda>0;\quad \min\{w_1^{\lambda},...,w_n^{\lambda}\}>0 \quad \mbox{ in } \Sigma_{\lambda}\quad  \}. $$
\begin{lem}\label{11nov17lem1}
$\bar \lambda=0$.
\end{lem}

\noindent{\bf Proof of Lemma \ref{11nov17lem1}:}
 If $\bar \lambda>0$, we first prove that $w_i^{\bar \lambda}>0$ in $\Sigma_{\bar \lambda}$ for all $i\in I$. Indeed, let $I_0=\{i\in I;\quad w_i^{\bar\lambda}\equiv 0\}$. If $I_0$ is not empty, the irreducibility of $A$ implies all $w_i^{\bar \lambda}\equiv 0$ in $\Sigma_{\bar \lambda}$. However, not all $\beta_i$ are $0$, so for some $i\in I$, we have
 $$\Delta w_i^{\bar \lambda}+\sum_{j\in I} a_{ij} |x|^{\beta_j}e^{\xi_j^{\bar \lambda}}w_j^{\bar \lambda}=-
 \sum_{j\in I} a_{ij}(|x^{\bar \lambda}|^{\beta_j}-|x|^{\beta_j})e^{u_j^{\bar\lambda}}<0. $$
 A contradiction.

 Next we derive a contradiction to the definition of $\bar\lambda$. Let $\lambda_k$ tend to  $\bar \lambda$ from the left. Thus $\lambda_k>0$ for all large $k$. We can assume that $\min_{i\in I} w_i^{\lambda_k}<0$ in $\Sigma_{\lambda_k}$ because otherwise, the strong maximum principle implies $w_i^{\lambda_k}>0$ in $\Sigma_{\lambda_k}$, a contradiction to the definition of $\bar \lambda$. Therefore, let $x_k$ be where the minimum of $\min_{i\in I}w_i^{\lambda}$ be attained and there is $i_k\in I$ such that $w_{i_k}^{\lambda_k}(x_k)=\min_{i\in I,x\in \Sigma_{\lambda_k}}w_i^{\lambda_k}<0$. By Lemma \ref{lesslem1}, $x_k\in B_R$ for some $R>0$ and all $k$. Along a subsequence $\{x_k\}$ converges to $\bar x\in \Sigma_{\bar \lambda}$ such that for some $i_0\in I$, $w_{i_0}^{\bar \lambda}(x_0)=0$. Since we have proved that $w_i^{\bar \lambda}>0$ for all $i\in I$ in $\Sigma_{\bar \lambda}$, $x_0\in T_{\bar \lambda}$. However, $\nabla w_{i_k}(x_k)=0$ leads to
 $\nabla w_{i_0}^{\bar \lambda}(x_0)=0$, a contradiction to the Hopf Lemma. Lemma \ref{11nov17lem1} is established. $\Box$

 \medskip

Thus we have proved $\bar \lambda=0$, which leads to
$$u_i(-x_1,x_2)\ge u_i(x_1,x_2), \quad \forall x_1\ge 0,\quad i\in I.  $$
Moving the plane from all possible directions we obtain the symmetry of $u_i$. Theorem \ref{beta0} is established. $\Box$

\medskip

\section{Uniqueness theorem on the linearized system}

In this section we prove Theorem \ref{nondegen}. The following lemma describes the projection of $u$ on $\sin k\theta$ and $\cos k\theta$.

\begin{lem}\label{highfre} Let $\phi_{i,k}(r)$ satisfy
\begin{equation}\label{11oct29e1}
\phi_{i,k}''+\frac 1r\phi_{i,k}'+\sum_{j\in I}a_{ij}r^{\beta_j}e^{u_j(r)}\phi_{j,k}-\frac{k^2}{r^2}\phi_{i,k}=0, \quad 0<r<\infty
\end{equation}
and
\begin{equation}\label{11oct29e2}
|\phi_{i,k}(r)|\le Cr^k(1+r)^{-2k}, \quad \forall r>0, \quad k\ge 1.
\end{equation}
If there exists $f=(f_1,..,f_n)$ such that
\begin{equation}\label{111027e1}
f_i''+\frac 1r f_i'-\frac{k^2}{r^2}f_i+\sum_{j=1}^n a_{ij}r^{\beta_j}e^{u_j}f_j<0, \quad r>0
\end{equation}
and
\begin{equation}\label{111027e2}
f_i(r)>0,\quad \forall r>0, \quad \lim_{r\to 0}f_i(r)/r^k=\infty, \quad \lim_{r\to \infty} f_i(r)r^k=\infty.
\end{equation}
Then $\phi_{ik}\equiv 0$.
\end{lem}

\noindent{\bf Proof of Lemma \ref{highfre}:}
We only need to show $\phi_{ik}\le 0$. Suppose this is not the case. Then because of the assumptions on the decay rates, without loss of generality we assume
\begin{equation}\label{11oct27e3}
w_1(r_0)=\frac{\phi_{1,k}(r_0)}{f_1(r_0)}=\max_{i,r}\frac{\phi_{i,k}(r)}{f_i(r)}>0.
\end{equation}
Note that the maximum can be attained because of the decay assumptions on $\phi_{i,k}$ and (\ref{111027e2}).
The equation for $w_1$, after simple derivation, is
\begin{eqnarray*}
w_1''+(\frac{2f_1'}{f_1}+\frac 1r)w_1'+w_1(\frac{f_1''}{f_1}+\frac 1r\frac{f_1'}{f_1}-\frac{k^2}{r^2}+a_{11}r^{\beta_1}e^{u_1})\\
+\sum_{j=2}^n a_{1j}r^{\beta_k}e^{u_j}\frac{\phi_{k,j}}{f_1}=0.
\end{eqnarray*}
Near $r_0$, $w_1(r)>0$. Thus in the neighborhood of $r_0$, using (\ref{111027e1}) we have
$$w_1''+(\frac{2f_1'}{f_1}+\frac 1r)w_1'>\sum_{j=2}^n a_{1j}r^{\beta_j}e^{u_j}\frac{f_jw_1-\phi_{k,j}}{f_1}. $$
The left hand side of the above is non-positive when evaluated at $r_0$, while the right hand side is non-negative. A contradiction.
Lemma \ref{highfre} is established. $\Box$

\bigskip

\noindent{\bf Proof of Theorem \ref{nondegen}:}

Let $f_i=-u_i'(r)$. Direct computation shows that
$$f_i''+\frac 1rf_i'-\frac{1}{r^2}f_i+\sum_j a_{ij}r^{\beta_j}e^{u_j}f_j=\sum_j a_{ij} \beta_j r^{\beta_j-1}e^{u_j}. $$
Since all $\beta_i<0$. $f=(f_1,...,f_n)$ satisfies (\ref{111027e1}) and (\ref{111027e2}) for all $k\ge 2$. Let $\phi^k=(\phi_1^k,...,\phi_n^k)$ be the radial part of the projection onto, say, $\sin k\theta$. Then $\phi^k$ satisfies (\ref{11oct29e1}). Since $\phi^k$ is bounded, it is easy to apply standard ODE theorem to obtain that (\ref{11oct29e2}) also holds.
Thus all the projections on $\sin k\theta$ and $\cos k\theta$ are all zero for $\beta\le 0$.

Finally we prove that for the projection on $\sin\theta$ or $\cos \theta$ is also zero. Let $\phi^1=(\phi_{1,1},..,\phi_{1,n})$ be the projection of $\phi$ on $\sin\theta$. Then
we have
$$\phi_{1,i}''+\frac 1r\phi_{1,i}'-\frac 1{r^2}\phi_{1,i}+\sum_j a_{ij}r^{\beta_j}e^{u_j}\phi_{1,j}=0. $$
Since $\phi^1$ is bounded, the standard ODE theory implies that $\phi_{1,i}$ behaves like $O(1/r)$ at infinity and like $O(r)$ near $0$. We shall use $f=(-u_1',...,-u_n')$ as the function to majorize $\phi^1$. To apply the same argument as in the proof of Lemma \ref{highfre}, The problem is the maximum may tend to $0$ or infinity. We first prove that this can not happen at $0$:
\begin{equation}\label{11oct29e4}
 \phi_i/f_i \mbox{ is strictly increasing near $0$ if $\phi_i$ is positive near $0$}.
 \end{equation}
 Clearly once (\ref{11oct29e4}) is proved, $\phi^1\equiv 0$, thus Theorem \ref{nondegen} would be established.

 Now we prove (\ref{11oct29e4}). Let $z_i=\phi_{1,i}/r$ and $F_i=f_i/r$. Direct computation yields
 $$z_i''+\frac 3rz_i'+\sum_j a_{ij}r^{\beta_j}e^{u_j}z_j=0, \quad r>0. $$
 and
 $$F_i''+\frac 3r F_i'+\sum_j a_{ij} r^{\beta_j}e^{u_j} F_j= \sum_j a_{ij} \beta_j r^{\beta_j-2} e^{u_j}. $$
 Since $\phi_{1,i}$ is positive near $0$, $z_i(0)>0$ (if $z_i(0)=0$, there is no need to consider this case, as the maximum can not tend to $0$).
 Easy to see that near $0$,
 $$z_i(r)=z_i(0)+\sum_j O(r^{\beta_j+2}) $$
 and
 $$F_i(r)=F_i(0)+\sum_j \frac{a_{ij}}{\beta_j+2}r^{\beta_j}e^{u_j(0)}+\sum_j O(r^{\beta_j+1}). $$
 Proving $\phi_{1,i}/f_i$ to be increasing near $0$ is the same as proving that $z_i/F_i$ is increasing near $0$. Since $\beta_i<0$, one immediately sees that $z_i/F_i$ is increasing near $0$.

 \medskip

 Next we prove that $z_i/F_i$ is decreasing if $z_i$ is positive at infinity. Assume $z_i(r)=q_i/r^2+O(r^3)$ at infinity. We have known that, for some $\delta_i>0$,
 $$u_i(r)=-m_i\log r+c_i +O(r^{-\delta_i})\quad r>1. $$
 Thus
 $$e^{u_i(r)}=e^{c_i}r^{-m_i}+O(r^{-m_i-\delta_i}), \quad r>>1. $$
 We obtain, by integration on the equation for $z_i$, that
 $$z_i'(r)=-\frac{2q_i}{r^3}+\sum_j a_{ij}\frac{e^{c_j}q_j}{m_j-\beta_j-2}r^{\beta_j-m_j-1}+O(r^{\beta_j-m_j-1-\delta_j}). $$
 $$z_i(r)=\frac{q_i}{r^2}-\sum_j a_{ij} \frac{e^{c_j}q_j}{(m_j-\beta_j-2)(m_j-\beta_j)}r^{\beta_j-m_j}+\sum_j O(r^{\beta_j-m_j-\delta_j}). $$

 Correspondingly to compute $F_i$, we use the equation for $u_i$ to obtain
 $$(ru_i'(r))'=-\sum_j a_{ij}r^{\beta_j+1}e^{u_j}=-\sum_j a_{ij}r^{1+\beta_j-m_j}e^{c_j}+O(r^{1+\beta_j-m_j-\delta_j}). $$
 Using $ru_i'(r)\to -m_i$ at infinity, we have
 $$ru_i'(r)=-m_i+\sum_j a_{ij}\frac{e^{c_j}}{m_j-\beta_j-2}r^{2+\beta_j-m_j}+O(r^{2+\beta_j-m_j-\delta_j}). $$
 Consequently
 $$F_i=-\frac{u_i'(r)}r =\frac{m_i}{r^2}-\sum_j a_{ij}\frac{e^{c_j}}{m_j-\beta_j-2}r^{\beta_j-m_j}+O(r^{\beta_j-m_j-\delta_j}). $$
 $$F_i'=-\frac{2m_i}{r^3}+\sum_j a_{ij}e^{c_j}\frac{m_j-\beta_j}{m_j-\beta_j-2}r^{\beta_j-m_j-1}+O(r^{\beta_j-m_j-1-\delta_j}). $$
 Our goal is
 $$z_i'F_i-z_iF_i'<0 \mbox{ near infinity } $$
 when $z_i$ is positive near infinity.
 Using the expressions above it is enough to show if the following is negative:
 \begin{equation}\label{12aug31e1}
 \sum_j a_{ij} e^{c_j}(-q_i+\frac{m_i}{m_j-\beta_j}q_j)r^{\beta_j-m_j-2}.
 \end{equation}
 By $q_i/m_i=\max_{j\in I} q_j/m_j$, $q_i>0$ and $\beta_i<0$ for all $i\in I$, we have (\ref{12aug31e1}). Therefore $z_i/F_i$ is decreasing near infinity. Theorem \ref{nondegen} is proved. $\Box$

\section{Appendix}

In this appendix, we prove the ODE lemmas stated in section two.

\noindent{\bf Proof of Lemma \ref{global1}:} The proof is standard ( for example, see \cite{chen-li-duke2}). We include it for the convenience of the reader. Let
\begin{equation}\label{11nov18e1}
w_i(x)=\int_{\mathbb R^2}(-\frac 1{2\pi}\log |x-y|+\frac 1{2\pi}\log (1+|y|))\sum_j a_{ij}|y|^{\beta_j}e^{u_j(y)}dy.
\end{equation}
Clearly $w_i$ is well defined and satisfies
$$-\Delta w_i(x)=\sum_j a_{ij} |x|^{\beta_j}e^{u_j(x)}, \quad \mathbb R^2$$
and
$$\Delta (u_i-w_i)=0, \quad \mathbb R^2. $$
By Lemma 4.1 in \cite{linzhang1} $u_i\le C$ on $\mathbb R^2$. Next we claim that
$$\lim_{|x|\to \infty} \frac{w_i(x)}{\log |x|}=-m_i. $$
To see the above, it is easy to obtain for $\epsilon>0$, there exist $R(\epsilon)>>1$ and $R_1>>R$ such that
for $|x|>R_1$
$$|\frac 1{2\pi}\int_{B_R}\frac{-\log |x-y|+\log (1+|y|)}{\log |x|}\sum_j a_{ij}|y|^{\beta_j}e^{u_j}-m_i|\le \epsilon .$$
Also
$$|\int_{\mathbb R^2\setminus B_R}\frac{-\log |x-y|+\log (1+|y|)}{\log |x|}\sum_j a_{ij}|y|^{\beta_j}e^{u_j}dy|\le \epsilon. $$
Thus $u_i-w_i\le C\log (1+|x|)$, which leads to
\begin{equation}\label{11nov18e2}
u_i=w_i+C_i
\end{equation} for some $C_i\in \mathbb R$.
Next we claim that
\begin{equation}\label{11sep14e1}
m_i-\beta_i>2 \quad \mbox{ for all } i\in I.
\end{equation}
 Indeed, if this is not the case, there exists $i_0\in I$ such that
$$m_{i_0}-\beta_{i_0}=\min\{m_1-\beta_1,..,m_n-\beta_n\}\le 2. $$
By (\ref{11nov18e1}) and (\ref{11nov18e2}) we have
$$u_{i_0}(x)=\frac{1}{2\pi}\int_{\mathbb R^2}(-\log |x-y|+\log (1+|y|))\sum_j a_{i_0 j}|y|^{\beta_j}e^{u_j(y)}dy-C. $$
Easy to check
$$-\log |x-y|+\log (1+|y|)\ge -\log (1+|x|), $$
thus
$$u_{i_0}(x)\ge -m_{i_0}\log (|x|+1)-C\ge -(2+\beta_{i_0})\log (|x|+1)-C$$
a contradiction to $\int_{\mathbb R^2}|x|^{\beta_{i_0}}e^{u_{i_0}(x)}<\infty$. (\ref{11sep14e1}) is established.
Now $u_i(x)$ can be written as
\begin{equation}\label{11nov10e1}
u_i(x)=-\frac 1{2\pi}\int_{\mathbb R^2}\log |x-y|\sum_j a_{ij}|y|^{\beta_j}e^{u_j(y)}dy+c_i.
\end{equation}
$c_i$ can be determined as in the statement. Finally we derive the error term $O(r^{-\delta})$. To see this we set
$$E_1=\{y;\quad |y|\le |x|/2. \}\quad E_2=\{y; \quad |y-x|\le \frac{|x|}2. \}\quad E_3=\mathbb R^2\setminus (E_1\cup E_2). $$
Using $e^{u_i(y)}=O(|y|^{-m_i})$ in (\ref{11nov10e1}) one obtains
$$-\frac 1{2\pi}\int_{E_1}\log |x-y|\sum_j a_{ij}|y|^{\beta_j}e^{u_j(y)}dy=-m_i\log |x|+O(|x|^{-\delta}). $$
Similarly by elementary estimates
$$\int_{E_2\cup E_3} \log |x-y|\sum_j a_{ij}|y|^{\beta_j}e^{u_j(y)}dy= O(|x|^{-\delta}). $$
The gradient estimate for $u_i$ is obtained by standard estimates.
Lemma \ref{global1} is established. $\Box$

\medskip

\noindent{\bf Proof of Lemma \ref{lineargrowth}:}
Let
$\psi(t)=(\psi_1(t),..., \psi_n(t))$ be defined as
$$\psi_i(t)=\phi_i(e^t),\quad i\in I. $$
Then $\psi$ satisfies
$$\psi_i''(t)+\sum_ja_{ij}e^{u_j(e^t)+(2+\beta_j)t}\psi_j(t)=0,\quad
-\infty<t<\infty,\quad i\in I. $$ Let
$\psi_{n+1}=\psi_1'$,..., $\psi_{2n}=\psi_n'$
and ${\mathbf F}=(\psi_1,..,\psi_{2n})^T$, then $\mathbf{F}$
satisfies
$$\mathbf{F}'=\mathbf{MF}$$
where $\displaystyle{\mathbf{M}=\left(\begin{array}{cc}
\mathbf{0}&\mathbf{I}\\
\mathbf{B}&\mathbf{0}
\end{array}
\right) }$. $\mathbf{B}$ is a $n\times n$ matrix with
$\mathbf{B}_{ij}=-a_{ij}e^{u_j(e^t)+(2+\beta_j)t}$. For $t>1$, the solution
for $\mathbf{F}$ is
\begin{equation}\label{jan22e1}
\mathbf{F}(t)=\lim_{N\to \infty}e^{\epsilon \mathbf{M}(t_N)}...e^{\epsilon \mathbf{M}(t_0)}\mathbf{F}(0).
\end{equation}
where $t_0,...,t_N$ satisfy $t_j=j\epsilon$, $j=0,..,N$, $\epsilon=t/N$.
 Since $u_i(e^t)+(2+\beta_i)t \sim (-m_i+2+\beta_i)t$ when $t$ is large and
$m_i>2+\beta_i$, we have $\|\mathbf{B}\|\sim
e^{-\delta t}$ for some $\delta>0$ and $t$ large. With this property we further have
\begin{equation}\label{jan22e2}
\|\mathbf{M}\|^k\le Ce^{-k\delta_1t},\quad k=2,3,... \quad t>0
\end{equation}
for some $\delta_1>0$.
Using (\ref{jan22e2}) in (\ref{jan22e1}) we have
$$\|\mathbf{F}(t)\|=O(t),\quad t>1. $$
Lemma \ref{lineargrowth} is established $\Box$

\medskip

\noindent{\bf Proof of Lemma \ref{ode1}:} If a solution $u=(u_1,...,u_n)$ exists, it would satisfy
$$u_i(r)=c_i(0)-\sum_{j=1}^n a_{ij}\int_0^r t^{\beta_j+1}(\log r-\log t)e^{u_j(t)}dt, \quad i=1,...,n. $$
We first prove the existence of a solution on $0<r<\delta$ for some small $\delta>0$ by iteration:
Let $u^{(0)}=(0,..,0)$ and
$$u_i^{(k+1)}(r)=c_i(0)-\sum_{j=1}^n a_{ij}\int_0^r t^{\beta_j+1}(\log r-\log t)e^{u_j^{(k)}(t)}dt, \quad i=1,...,n. $$
For $\delta>0$ small and $r\in (0,\delta)$, since $\beta_i+1>-1$, it is easy to see that such a sequence converges. Therefore the existence of a solution for over $(0,\delta)$ is proved.
The existence for $r\in (\delta, \infty)$ clearly holds because of the right hand side is a Lipschitz function of $u$. The proof of the
uniqueness of the solution is the same as that in Lemma \ref{linearlem2} later in the section. Lemma \ref{ode1} is established. $\Box$

\medskip

\noindent{\bf Proof of Lemma \ref{ode2}:} By Lemma \ref{ode1} a solution to (\ref{11sep9e1}) exists for $r>0$. We just need to show that
$$\int_0^{\infty}e^{u_i(r)}r^{\beta_i+1}dr<\infty,\quad \forall i\in I. $$
Let $v_i(t)=u_i(e^t)+(2+\beta_i)t$ ($i\in I$), then $v=(v_1,...,v_n)$ satisfies
$$v_i''(t)+\sum_j a_{ij} e^{v_j(t)}=0, \quad -\infty<t<\infty, \quad i\in I. $$
From the equation for $u$ we have
$$ru_i'(r)=-\int_0^r \sum_j a_{ij} e^{u_j(s)}s^{\beta_j+1}ds<0,\quad r>0, \quad i\in I. $$
The last inequality is strict because $a_{ij}\ge 0$ and not all equal to $0$. Consequently $v_i'(t)<2+\beta_i$ for $t\in \mathbb R$. Fix $t_0\in \mathbb R$ we have, for $t>t_0$,
$$v_i'(t)=v_i'(t_0)-\int_{t_0}^t \sum_j a_{ij} e^{v_j(s)}ds\le v_i'(t_0)-a_{ii}\int_{t_0}^t e^{v_i(s)}ds, \quad i\in I. $$
Since $a_{ii}>0$ there exists $t>t_0$ such that $v_i'(t)<0$. Choose $v_i'(t_1)=-\delta<0$ for some $\delta>0$, then we see
$$v_i(t)\le v_i(t_1)-\delta(t-t_1), \quad t>t_1 $$
which is equivalent to $u_i(r)<(-2-\beta_i-\delta)\log r+C$ for $r>e^{t_1}$. Therefore $\int_0^{\infty}e^{u_i(r)}r^{\beta_i+1}dr<\infty. $ Lemma \ref{ode2} is established. $\Box$

\medskip

\noindent{\bf Proof of Lemma \ref{linearlem2}:}
The proof is standard, we include it for the convenience of the reader. Clearly we only need to show that $\phi_i\equiv 0$ in $(0,\delta)$ for $\delta>0$ small. Write $\phi_i(r)$ as
\begin{eqnarray*}
\phi_i(r) &=&-\int_0^r\frac 1s \int_0^s (\sum_{j\in I} a_{ij}t^{\beta_j+1}e^{u_j(t)}\phi_j(t)dt)ds\\
&=&-\int_0^r\sum_{j\in I}a_{ij}t^{\beta_j+1}e^{u_j(t)}\phi_j(t)(\log r-\log t)dt
\end{eqnarray*}
Let $\alpha=\min\{\beta_1,...,\beta_n\}+1$, since all $\beta_i>-2$ we have $\alpha-\epsilon>-1$ for some $\epsilon>0$ small. For the $\epsilon$ we choose $\delta>0$ small so that $\log r-\log t<t^{-\epsilon}$ for $r<\delta$ and $t\le r$. Thus
$$|\phi_i(r)|\le C\int_0^r t^{\alpha-\epsilon}\sum_{j}|\phi_j(t)|dt, \quad r<\delta $$
for some $C$.
Let $\phi(r)=\sum_{i\in I}|\phi_i(r)|$ and $F(r)=\int_0^r t^{\alpha-\epsilon}|\phi(t)|dt$,
then
$$F'(r)-Cr^{\alpha-\epsilon}F(r)\le 0,\quad F\ge 0, \quad F(0)=0. $$
Since $\alpha-\epsilon>-1$,
$F\equiv 0$. Lemma \ref{linearlem2} is established. $\Box$


\begin{thebibliography}{99}

\bibitem{bartolucci1} D. Bartolucci, Compactness result for periodic multivortices in the electroweak theory, Nonlinear
Anal., 53 (2003), 277–297.

\bibitem{bartolucci2}
D. Bartolucci and G. Tarantello, The Liouville equation with singular data: A concentration-
compactness principle via a local representation formula, J. Differential Equations, 185
(2002), 161–180.

\bibitem{bartolucci3} D. Bartolucci and G. Tarantello, Liouville type equations with singular data and their appli-
cations to periodic multivortices for the electroweak theory, Comm. Math. Phys., 229 (2002),
3–47.

\bibitem{bclt}  D. Bartolucci, C. C. Chen, C. S. Lin and G. Tarantello, Profile of blow-up solutions to mean
field equations with singular data, Comm. Partial Differential Equations, 29 (2004), 1241–
1265.


\bibitem{bennet} W. H. Bennet, Magnetically self-focusing streams, Phys. Rev. 45 (1934),
890–897.
\bibitem{BM} H. Brezis and F. Merle, Uniform estimates and blow-up behavior for solutions of $-\Delta u=V(x)e^u$ in two dimensions, Comm. Partial Differential Equation, 16 (1991), 1223–1253.

\bibitem{caff-yang} L. Caffarelli and Y. Yang, Vortex condensation in the Chern-Simons Higgs model: an exis-
tence theorem, Comm. Math. Phys., 168 (1995), 321–336.

\bibitem{caglioti1} E. Caglioti, P. L. Lions, C. Marchioro and M. Pulvirenti, A special class of stationary flows for
two-dimensional Euler equations: A statistical mechanics description, Comm. Math. Phys.,
143 (1992), 501–525.

\bibitem{caglioti} E. Caglioti, P. L. Lions, C. Marchioro and M. Pulvirenti, A special class of stationary flows
for two-dimensional Euler equations: A statistical mechanics description, part II, Comm.
Math. Phys., 174 (1995), 229–260.

\bibitem{chanillo1} S. Chanillo and M. Kiessling, Rotational symmetry of solutions of some nonlinear problems
in statistical mechanics and in geometry, Comm. Math. Phys., 160 (1994), 217–238.

\bibitem{chanillo2}S. Chanillo, M. K-H Kiessling,
Conformally invariant systems of nonlinear PDE of Liouville type.
Geom. Funct. Anal. 5 (1995), no. 6, 924--947.

\bibitem{chang1} S. Y. Chang, M. Gursky and P. Yang, The scalar curvature equation on 2- and 3-spheres,
Calc. Var. and PDE, 1 (1993), 205–229.

\bibitem{chang2} S. Y. Chang and P. Yang, Prescribing Gaussian curvatuare on S2, Acta Math., 159 (1987),
215–259.

\bibitem{chenlin0} C. C. Chen and C. S. Lin, Estimate of the conformal scalar curvature equation via the method
of moving planes. II, J. Differential Geom., 49 (1998), 115–178.

\bibitem{chenlin1} C. C. Chen, C. S. Lin, Sharp estimates for solutions of multi-bubbles in compact Riemann surfaces,
Comm. Pure Appl. Math. 55 (2002), no. 6, 728--771.

\bibitem{chenlin2} C. C. Chen, C. S. Lin, Topological degree for a mean field equation on Riemann surfaces,
Comm. Pure Appl. Math. 56 (2003),  no. 12, 1667--1727.

\bibitem{chenlin3} C. C. Chen, C. S. Lin, Mean field equations of Liouville type with singular data: Sharper estimates,
Discrete and continuous dynamic systems, Vol 28, No 3, 2010, 1237-1272.

\bibitem{chenlinwang} C. C. Chen, C. S. Lin and G. Wang, Concentration phenomena of two-vortex solutions in a
Chern-Simons model, Ann. Sc. Norm. Super. Pisa Cl. Sci. (5), 3 (2004), 367–397.

\bibitem{cheng-li-duke1} W. Chen and C. Li, Classification of solutions of some nonlinear elliptic equations, Duke
Math. J., 63 (1991), 615–622.

\bibitem{chen-li-duke2} W. X. Chen, C. M. Li, Qualitative properties of solutions to some nonlinear elliptic equations in $\mathbb R^2$.
Duke Math. J. 71 (1993), no. 2, 427–439.

\bibitem{childress} S. Childress and J. K. Percus, Nonlinear aspects of Chemotaxis, Math. Biosci. 56 (1981),
217–237.

\bibitem{CSW} M. Chipot, I. Shafrir, G. Wolansky, On the solutions of Liouville systems.  J. Differential Equations  140  (1997),  no. 1, 59--105.

\bibitem{CSW1} M. Chipot, I. Shafrir, G. Wolansky, Erratum: ``On the solutions of Liouville systems'' [J. Differential Equations 140 (1997), no. 1, 59--105; MR1473855 (98j:35053)].  J. Differential Equations  178  (2002),  no. 2, 630.

\bibitem{debye} P. Debye and E. Huckel, Zur Theorie der Electrolyte, Phys. Zft 24 (1923), 305–325.

\bibitem{dunne} G. Dunne, Self-dual Chern-Simons Theories, Lecture
Notes in Physics, vol. m36, Berlin: Springer-Verlag, 1995.

\bibitem{hong} J. Hong, Y. Kim and P. Y. Pac, Multivortex solutions of the Abelian Chern-Simons-Higgs
theory, Phys. Rev. Letter, 64 (1990), 2230–2233.

\bibitem{jackiw} R. Jackiw and E. J. Weinberg, Selfdual Chern Simons vortices, Phys. Rev. Lett., 64 (1990),
2234–2237.

\bibitem{jost-lin-wang} J. Jost, C. S. Lin and G. Wang, Analytic aspects of the Toda system. II. Bubbling behavior
and existence of solutions, Comm. Pure Appl. Math., 59 (2006), 526–558.

\bibitem{jw1} J. Jost, G. Wang, Classification of solutions of a Toda system in ${\mathbb R}\sp 2$.
 Int. Math. Res. Not. 2002, no. 6, 277--290.

\bibitem{jw2} J. Jost, G. Wang, Analytic aspects of the Toda
system. I. A Moser-Trudinger inequality. Comm. Pure Appl. Math. 54
(2001), no. 11, 1289--1319.

\bibitem{kw} J. Kazdan and F. Warner, Curvature functions for compact 2-manifolds, Ann. of Math., 99
(1974), 14–47.

\bibitem{keller} E. F. Keller and L. A. Segel, Traveling bands of Chemotactic Bacteria: A theoretical analysis,
J. Theor. Biol. 30 (1971), 235–248.

\bibitem{kiessling2}M. K.-H. Kiessling and J. L. Lebowitz, Dissipative stationary Plasmas: Kinetic Modeling Bennet
Pinch, and generalizations, Phys. Plasmas 1 (1994), 1841–1849.

\bibitem{licmp} Y. Y. Li, Harnack type inequality: the method of moving planes,
Comm. Math. Phys. 200 (1999), no. 2, 421--444.

\bibitem{lin1} C. S. Lin, An expository survey on the recent development of mean field equations, Discrete
Contin. Dyn. Syst., 19 (2007), 387–410.

\bibitem{lin2} C. S. Lin and M. Lucia, Uniqueness of solutions for a mean field equation on torus, J.
Differential Equations, 229 (2006), 172–185.

\bibitem{lin3} C. S. Lin and M. Lucia, One-dimensional symmetry of periodic minimizers for a mean field
equation, Ann. Sc. Norm. Super. Pisa Cl. Sci. (5), 6 (2007), 269–290.

\bibitem{lwy} C. S. Lin, J. C. Wei, D. Ye, Classifcation and nondegeneracy of $SU(n+1)$ Toda system. preprint.

\bibitem{lwz} C. S. Lin, J. C. Wei, C. Zhao, Sharp estimates for fully bubbling solutions of a $SU(3)$ Toda system. Preprint.

\bibitem{linzhang1} C. S. Lin, L. Zhang, Profile of bubbling solutions to a Liouville system.  (English, French summary) Ann. Inst. H. Poincaré Anal. Non Linéaire 27 (2010), no. 1, 117--143.

\bibitem{linzhang2} C. S. Lin, L. Zhang,   A topological degree counting for some Liouville systems of mean field equations, Comm. Pure Appl. Math. volume 64, Issue 4, pages 556–-590, April 2011.

\bibitem{mock} M. S. Mock, Asymptotic behavior of solutions of transport equations for semiconductor devices,
J. Math. Anal. Appl. 49 (1975), 215–225.

\bibitem{nolasco1} M. Nolasco and G. Tarantello, On a sharp type inequality on two dimensional compact man-
ifolds, Arch. Rational Mech. Anal., 145 (1998), 161–195.

\bibitem{nolasco2} M. Nolasco and G. Tarantello, Double vortex condensates in the Chern-Simons-Higgs theory,
Calc. Var. and PDE, 9 (1999), 31–94.

\bibitem{nolasco3} M. Nolasco and G. Tarantello, Vortex condensates for the SU(3) Chern-Simons theory,
Comm. Math. Phys., 213 (2000), 599–639.

\bibitem{prajapat} J. Prajapat and G. Tarantello, On a class of elliptic problems in R2: Symmetry and unique-
ness results, Proc. Roy. Soc. Edinburgh Sect. A, 131 (2001), 967–985.

\bibitem{ricciardi} T. Ricciardi and G. Tarantello, Vortices in the Maxwell-Chern-Simons theory, Comm. Pure
Appl. Math., 53 (2000), 811–851.

\bibitem{rubinstein} I. Rubinstein, Electro diffusion of Ions,
SIAM, Stud. Appl. Math. 11 (1990).

\bibitem{spruck-yang1} J. Spruck and Y. Yang, Topological solutions in the self-dual Chern-Simons theory: Existence
and approximation, Ann. Inst. H. Poinc´are Anal. Non Lin´eaire, 12 (1995), 75–97.

\bibitem{spruck-yang2} J. Spruck and Y. Yang, On multivortices in the electroweak theory I: Existence of periodic
solutions, Comm. Math. Phys., 144 (1992), 1–16.

\bibitem{spruck-yang3} J. Spruck and Y. Yang, On multivortices in the electroweak theory II: Existence of Bogo-
mol'nyi solutions in R2, Comm. Math. Phys., 144 (1992), 215–234.

\bibitem{tarantello1} G. Tarantello, A Harnack inequality for Liouville-type equations with singular sources, Indiana
Univ. Math. J., 54 (2005), 599–615.

\bibitem{tarantello2} G. Tarantello, A quantization property for blow up solutions of singular Liouville-type equa-
tions, J. Funct. Anal., 219 (2005), 368–399.

\bibitem{tarantello3} G. Tarantello, Analytical aspects of Liouville-type equations with singular sources, Stationary
partial differential equations, Vol. I, 491–592, Handb. Differ. Equ., North-Holland, Amsterdam,
2004.

\bibitem{tarantello4} G. Tarantello, Selfdual gauge field vortice. An analytical approach, in "Contributions in Nonlinear
Differential Equations and their Applications," 72, Birkh¨auser, Boston, 2008.

\bibitem{yang} Y. Yang, "Solitons in Field Theory and Nonlinear Analysis," Springer-Verlag, New York, 2001.

\bibitem{zhangcmp} L. Zhang, Blow up solutions of some nonlinear elliptic equations involving exponential non-
linearities, Comm. Math. Phys., 268 (2006), 105–133.

\bibitem{zhangccm} L. Zhang, Asymptotic behavior of blowup solutions for elliptic equations with exponential
nonlinearity and singular data, Communications in Contemporary mathematics, vol 11, issue 3, June 2009, 395-411.



\end{thebibliography}
\end{document}